# Mathematical and computational approaches for stochastic control of river environment and ecology: from fisheries viewpoint

Hidekazu Yoshioka[1],*

**Abstract**  We present a modern stochastic control framework for dynamic optimization of river environment and ecology. We focus on a fisheries problem in Japan and show several examples of simplified optimal control problems of stochastic differential equations modeling fishery resource dynamics, reservoir water balance dynamics, benthic algae dynamics, and sediment storage dynamics. These problems concern different phenomena with each other, but they all reduce to solving degenerate parabolic or elliptic equations. Optimal controls and value functions of these problems are computed using finite difference schemes. Finally, we present a higher-dimensional problem of controlling a dam-reservoir system using a semi-Lagrangian discretization on sparse grids. Our contribution shows the state-of-art of modeling, analysis, and computation of stochastic control in environmental engineering and science, and related research areas.

## 1. Introduction

Rivers are a part of hydrological cycles as well as a part of human lives. Flowing waters in rivers are stored by dam-reservoir systems and utilized as primary water resources for drinking, irrigation, hydropower generation, and so on [58]. Operation of dam-reservoir systems should harmonize the water use with river water environment and ecosystems. Regulated flows released from a dam alter its downstream flow regimes and often threaten riparian habitats and aquatic species [27, 70].

[1] H. Yoshioka (The corresponding author)
Graduate School of Natural Science and Technology, Shimane University and Fisheries Ecosystem Project Center, Shimane University
Nishikawatsu-cho 1060, Matsue, 690-8504, Japan
yoshih@life.shimane-u.ac.jp

Therefore, exploring a unified framework to balance ecological and human dimensions is of high importance in river environmental management.

Stochastic optimal control as a branch of modern mathematical sciences plays a central role in analysis and management of noise-driven dynamical systems [42]. Noises in the context of environmental and ecological management arise from stochastic river water flows [56] and highly nonlinear and possibly unresolved biological phenomena such as biological growth phenomena [80]. Stochastic differential equations (SDEs) [42], which are formally seen as ordinary differential equations (ODEs) driven by noises, serve as an efficient mathematical tool for modeling and controlling noisy system dynamics.

Problems related to river environmental and ecological dynamics are not the exception that the stochastic control applies. However, such an outlook has not been paid much attention before the author and his co-workers started modeling, analysis, and computation of inland fishery resource management in rivers in Japan. Fishery resources management problems in seas have conventionally been studied as both stochastic and deterministic optimal control problems [20, 50, 32, 33], possibly because of their huge impacts on food and economy worldwide. By contrast, the problems in inland waters usually have smaller impacts; nevertheless, they have been serving as unique elements to sustain local ecosystems, societies, and sometimes ecological education [76]. In addition, as we will demonstrate in this chapter, there exist many interesting management problems specific to inland fishery resources. These problems have been emerging as a new application of the stochastic control to engineering problems.

The objective of this chapter is to present a unified mathematical framework and specific applications of the stochastic control to river environmental and ecological problems with a focus on inland fisheries management. Problems we focus on concern management of the diadromous fish *Plecoglossus altivelis altivelis* (*P. altivelis*, Ayu) in Japan as one of the most common inland fishery resources in the country [1, 37]. The unique life history of *P. altivelis* is explained later, but what is important here is that managing the fish requires considering not only its life history, but also surrounding environmental conditions from multiple sides. This motivates us to separately study sub-control problems, such as fish growth [77], benthic algae management [73], sediment storage management [74], and dam-reservoir system management [86], which have different characteristics with each other but can be handled by a unified framework based on dynamic programming and viscosity solutions: appropriate weak solutions to degenerate elliptic equations [3, 16].

We show that each control problem reduces to solving a corresponding degenerate elliptic (or parabolic, hyperbolic) equation often called Hamilton-Jacobi-Bellman (HJB) equation. This is carried out either analytically or numerically, but usually the latter is employed in applications because of the nonlinearity and nonlocality of the HJB equations [7, 22, 24, 34, 51]. We approach the specific problem numerically using finite difference and semi-Lagrangian schemes.

We also consider a coupled higher-dimensional problem where resource and environmental dynamics should be managed concurrently. Conventional numerical schemes encounter a huge computational cost when there are more than three to four state variables. This issue is called the curse of dimensionality. Computational costs of the conventional numerical methods for stochastic control problems increase exponentially when the total number of state variables increases linearly. We use a semi-Lagrangian scheme [8] on sparse grids [10] to alleviate this issue.

This chapter is organized as follows. The conventional framework of stochastic control is briefly reviewed in the second section. The specific problems, which are based on our studies but are slightly reformulated here for the sake of consistency, are separately studied in the third section. Their applicability and limitation are discussed as well. A coupled problem is analyzed in the fourth section. Our summary and future perspectives are presented in the last section.

## 2. Stochastic control

### 2.1 Stochastic differential equation

We briefly and formally present a basic framework of the stochastic control. Interested readers should read the textbooks and reviews of stochastic control and its applications [12, 42, 71, 72]. Throughout this chapter, each problem is considered based on a standard complete probability space [42]. The explanation in this section is formal, and coefficients and parameters appearing in the third section will be specified in each control problem.

We consider continuous-time dynamics. The time is denoted as $t \geq 0$. The total number of state variables is denoted as $M \in \mathbb{N}$. The state variables are assumed to be càdlàg and are represented as a vector $\mathbf{X}_t = \left[ X_{i,t} \right]_{1 \leq i \leq M}$. Its range is denoted as $\Omega$, which is bounded or unbounded, and is problem-dependent. We assume that the process $\mathbf{X} = \left( \mathbf{X}_t \right)_{t \geq 0}$ is a jump-diffusion process governed by a system of SDEs driven by compound Poisson jumps. From a mathematical side, it is more convenient to consider a generic Lévy process as a driving noise process [42]. However, problems considered in this chapter requires only compound Poisson processes (and possibly Brownian motions if necessary), which are assumed to be mutually independent with each other.

We only use jump processes in the mathematical modeling discussed later but incorporating diffusive dynamics does not encounter difficulties in most cases. Both Brownian and jump noises are considered in this sub-section for the sake of explanation. Recall that the continuous-time Markov chains are represented using appropriate Poisson processes [71].

The $N$-dimensional standard Brownian motion is denoted as $\mathbf{B}_t = \left[B_{i,t}\right]_{1 \leq i \leq N}$ at time $t$. In addition, the $K$-dimensional compound Poisson process is denoted as $\mathbf{P}_t = \left[P_{i,t}\right]_{1 \leq i \leq K}$ at time $t$. Each $P_{i,t}$ is mutually independent with each other. For each $1 \leq i \leq K$, the jump intensity of $P_{i,t}$ is denoted as $\lambda_i > 0$, and the probability density function of the jump size of $P_{i,t}$ as $p_i$. We assume $\mathbf{X}_t \in \Omega$ a.s. for $t \geq 0$. The control process $\mathbf{u} = (\mathbf{u}_t)_{t \geq 0}$ is assumed to have a compact range $U$ and progressively measurable with respect to a natural filtration generated by the processes $(\mathbf{B}_t)_{t \geq 0}$ and $(\mathbf{P}_t)_{t \geq 0}$. These assumptions are rather standard. The admissible set of $\mathbf{u}$ is denoted as $\mathcal{U}$. The system of Itô's SDEs governing $\mathbf{X}$ is set as

$$\mathrm{d}\mathbf{X}_t = b(t, \mathbf{X}_t, \mathbf{u}_t)\mathrm{d}t + \sigma(t, \mathbf{X}_t)\mathrm{d}\mathbf{B}_t + l(t, \mathbf{X}_{t-}, \mathbf{u}_{t-})\mathrm{d}\mathbf{P}_t, \ t > 0 \quad (1)$$

subject to an initial condition $\mathbf{X}_0 \in \Omega$. Here, $\mathbf{X}_{t-}$ is the left limit of $\mathbf{X}_t$ at time $t$, and each coefficients $b, \sigma, l$ are assumed to have suitable dimensions and chosen so that, for each control $\mathbf{u}$, the system (1) has a unique path-wise solution. Sufficient conditions for the unique existence are found in textbooks such as [12, Chapter 4]. We assume that only $b$ and $l$ are modulated by $\mathbf{u}$ to simplify the explanation.

## 2.2 Performance index and value function

The performance index is a functional of the time $t$, the process $\mathbf{X}$, and the control $\mathbf{u} \in \mathcal{U}$, to be minimized with respect to $\mathbf{u}$ by a decision-maker: the controller of the dynamics (1). The expectation conditioned on $\mathbf{X}_t = \mathbf{x}$ is denoted as $\mathbb{E}^{t,x}$. The terminal time is denoted as $T > 0$, which is possibly unbounded. The performance index $\phi = \phi(t, \mathbf{X}, \mathbf{u})$ conditioned on $\mathbf{X}_t = \mathbf{x}$ is set as

$$\phi(t, \mathbf{X}, \mathbf{u}) = \mathbb{E}^{t,x}\left[\int_t^T f(s, \mathbf{X}_s, \mathbf{u}_s) e^{-\delta(s-t)} \mathrm{d}s + g(T, \mathbf{X}_T) e^{-\delta(T-t)}\right]. \quad (2)$$

Here, $\delta \geq 0$ is the discount rate representing myopia of the decision-maker; a larger $\delta$ means that he/she is more myopic (e.g., Bian et al. [5]). The coefficients $f, g$ are sufficiently regular and bounded in $[0,T] \times \Omega \times U$ and $\Omega \times U$, respectively. The first and second terms of (2) mean the cumulative utility/disutility and the penalty incurred at $T > 0$. For a problem with an unbounded $T$ where $b, \sigma, l, f$ are time-independent and $g \equiv 0$, $\phi$ is alternatively set as

$$\phi(\mathbf{X}, \mathbf{u}) = \mathbb{E}^x\left[\int_0^\infty f(\mathbf{X}_s, \mathbf{u}_s) e^{-\delta s} \mathrm{d}s\right], \quad (3)$$

provided that the right-hand side exists, where $\mathbb{E}^x$ represents $\mathbb{E}^{0,x}$. Notice that in both the finite-horizon and infinite-horizon cases, we assume $\mathbf{X}_t \in \Omega$ a.s. for $t \geq 0$ so that the problem is not terminated in $[0,T]$.

The value function is the minimized $\phi$ with respect to $\mathbf{u} \in \mathcal{U}$:

$$\Phi(t,\mathbf{X}) = \inf_{\mathbf{u} \in \mathcal{U}} \phi(t,\mathbf{X},\mathbf{u}) \text{ for } T < +\infty \ (\Phi(\mathbf{X}) = \inf_{\mathbf{u} \in \mathcal{U}} \phi(\mathbf{X},\mathbf{u}) \text{ for } T = +\infty). \quad (4)$$

A minimizer in (4), if it exists, is referred to as an optimal control and is denoted as $\mathbf{u}^*$. Without significant loss of generality, we only consider Markov controls formally represented in a feedback form $\mathbf{u}_t^* = \mathbf{u}^*(t,\mathbf{X}_t)$ in the finite-horizon case and $\mathbf{u}_t^* = \mathbf{u}^*(\mathbf{X}_t)$ in the infinite-horizon case. The Markov control assumption is not restrictive in many applications (e.g., Øksendal and Sulem [42]).

### 2.3 HJB equation

The HJB equation is a nonlinear and possibly nonlocal degenerate parabolic differential equation formally derived as a governing equation of the value function $\Phi$. Based on a dynamic programming principle (e.g., Øksendal and Sulem [42], Pham [45], Touzi [62]), the HJB equation in the finite-horizon case becomes

$$\begin{aligned}
-\frac{\partial \Phi}{\partial t} + \delta \Phi - \frac{1}{2} \sigma_{ik} \sigma_{kj} \frac{\partial \Phi}{\partial x_i \partial x_j} & \\
-\inf_{\mathbf{u} \in U} \left\{ b_i \frac{\partial \Phi}{\partial x_i} + f - \lambda_i \int \left\{ \Phi - \Phi(t, \mathbf{x}(1 + l(t,\mathbf{u},\mathbf{x}))) \right\} p_i(z_i) \mathrm{d}z_i \right\} &= 0
\end{aligned} \text{ in } [0,T] \times \Omega \quad (5)$$

subject to the terminal condition $\Phi = g$ in $\Omega$ at $t = T$, where the Einstein's convention has been used in (5). Some boundary conditions may be prescribed along the boundary of $\Omega$. This is not explicitly considered here but will be prescribed in each problem below when necessary. The coefficients of the dynamics (1) and the performance index (2) are inherited in this HJB equation. The equation in the infinite-horizon case is derived by simply omitting the first term of (5) and the terminal condition, and setting the domain of the equation as $\Omega$.

### 2.4 Remarks

We close this section with several remarks on HJB equations. A candidate of the optimal control $\mathbf{u}_t^*$ is formally derived as

$$\mathbf{u}^*(t,\mathbf{x}) = \arg\min_{\mathbf{u} \in U} \left\{ b_i(t,\mathbf{x},\mathbf{u}) \frac{\partial \Phi}{\partial x_i}(t,\mathbf{x},\mathbf{u}) + f(t,\mathbf{x},\mathbf{u}) \right\}, \quad (6)$$

implying, at least formally, that solving the stochastic control problem ultimately reduces to finding the value function $\Phi$ by solving the HJB equation (5). This point is an advantage of the stochastic control model because we can obtain the optimal control $\mathbf{u}^*$ for all the possible states if we can solve the HJB equation once, depending on the performance index $\phi$ that can be determined flexibly. A disadvantage is that solving an HJB equation is not always easy, but the equation can be resolved numerically for moderately small systems having two to three state variables.

Key mathematical problems in the stochastic control modeling are existence, uniqueness, and regularity of solutions to the HJB equation (5). If the coefficients of the problem satisfy certain boundedness and regularity assumptions, then an HJB equation would admit a unique classical solution satisfying the equation pointwise (e.g., Bian et al. [6], Pham [44]). However, this is not the case in general. We often must seek for solutions in a weaker sense because of the loss of regularity of solutions where the regularity of some of the coefficients in the problem drop.

The most plausible candidate of weak solutions can be viscosity solutions [16], where one-sided HJB equations (Formally, " = " in (5) is replaced by " ≤ " or " ≥ ") are satisfied by appropriate sufficiently smooth test functions. The most useful point in relying on the concept of viscosity solutions is that solutions to an HJB equation need not be continuously differentiable, and even not required to be continuous in some cases (e.g., Touzi [62]). Notice that a classical solution is a viscosity solution (the converse statement is false). In this view, it is natural to analyze HJB equations from a viscosity viewpoint. The definition of viscosity solutions is not presented here but will be found in the references of the separate problems below. The existence and uniqueness of solutions to HJB equations are in general proven with the help of a comparison argument stating that a viscosity sub-solution is always not greater than a viscosity super-solution [16].

The control variable employed in the explanation above was assumed to be continuous-time; however, there exist other classes of control variables of importance in both theory and applications. Such examples include impulse controls and singular controls [42, Chapters 8 and 9] and partial observation controls [46]. In these problems, the basic strategy is essentially the same: set dynamics, set a performance index, derive an optimality equation, and finally solve it to find an optimal control. Their optimality equation is more complicated and more difficult to handle. We discuss one related simple but a delicate example at the end of Section 3.

In most cases, HJB equations are not solvable analytically, and must be solved numerically. Convergent numerical schemes, such as finite difference and semi-Lagrangian schemes, can be developed based on the monotonicity, stability, and consistency requirements [4]. Unfortunately, these schemes do not always exhibit satisfactory accuracy in applications. There have been recent progresses on high-resolution schemes for computing viscosity solutions [23, 47, 49]. Although not discussed here, variational characterizations of HJB equations are also possible [31], with which weak differentiability of solutions can be analyzed.

## 3. Specific problems

### *3.1 "Non-renewable" fishery resource management*

In this section, we present four problems related to inland fisheries and river environmental management. Each problem is formulated under a simplified setting to clearly present the structure of the problem, value function, and optimal control. Different numerical schemes are used in different problems, so that readers can understand that many schemes are available for computing HJB equations.

The first problem is planning a harvesting policy of a "non-renewable" fishery resource in a finite horizon. The fish considered here is *P. altivelis*: one of the major inland fishery resources in Japan as a core of the aquatic ecosystem as well as an indispensable element shaping regional culture [77]. Harvesting the fish is for recreational purposes in the country in most cases. A recent survey revealed that the fish is the third most popular inland fishery resource in Japan [39-40]. The fish is diadromous (especially, born in the sea and grow up in an upstream river) having a unique one-year life history. The explanation below is based on Yoshioka et al. [83]. See, also Murase et al. [38]. The life history of the fish is briefly explained as follows. Adult fishes grow up in mid-stream of a river during spring to the coming summer by feeding rock-attached algae like diatoms and migrate together toward the downstream river reach in the coming autumn to spawn. After the spawning, adult fishes die. The hatched larvae drift toward the sea and grow up by feeding mainly on plankton in the sea until coming spring, at which mass migrations of the fish population from the sea to rivers occur.

Conventionally, the harvesting season of the fish in a river in Japan is the early summer (June to July) to the late coming autumn (October to November) with slight seasonal differences among different rivers in the country. The fish as a fishery resource is renewable unless its life cycles are not terminated due to an extinction. However, it can be seen as a non-renewable fishery resource from the spring to the coming autumn during which the fish does not spawn and their population monotonically decreases by harvesting, natural deaths, and predation by the waterfowls such as the great cormorant *Phalacrocorax carbo* [69].

To achieve a sustainable fisheries management of the fish, one has not to harvest the fish too much, while the predation from the other species like the waterfowls should be mitigated if possible, in order to increase the population that can be potentially harvested. The decision-maker, which is a local fishery cooperative and its union members, should therefore concurrently consider both harvesting and protection of the fish.

The formulation below follows Yoshioka and Yaegashi [77, 79] but is modified so that the model becomes simpler but still non-trivial. The problem here considers harvesting a non-renewable (lumped) population in a finite interval $[0,T]$ with $T > 0$. The population at time $t$ is denoted as $X_t$. Physically, the population is the total number of individual fishes, and should be an integer variable. However, we regard it as a real variable assuming that the population is sufficiently large. This theoretical assumption allows us to formulate the mathematical modeling here in a more tractable manner. We assume that the population does not increase during $[0,T]$. We also assume that harvesting a larger individual fish is more profitable and contributes to a larger utility of the decision-maker. In this view, the body weight of the fish should also be considered. Biologically, it is quite natural to consider that individual fishes have different growth curves. Such phenomena can be described with the individual-based models [14] but can be more complicated than the model presented below. Therefore, we assume that the decision-maker knows the mean growth curve of the individual fishes, and set it as a Lipschitz continuous function $U_t$ at time $t$. This is assumed to be a given positive and increasing variable. Logistic-like, sigmoid-shape growth curves [81] can be used to represent $U_t$. The (mean) total biomass at $t$ is $X_t U_t$.

We assume that the population decreases due to natural death, predation, or harvesting. The natural deaths are assumed to be due to both time-continuous and catastrophic events, where the latter is described by a compound Poisson process [52] $P = (P_t)_{t \geq 0}$. The natural (and gradual) mortality rate is denoted as $R > 0$, the predation pressure as $p > 0$, the harvesting pressure as $h = (h_t)_{t \geq 0}$, and the jump intensity of the catastrophic deaths as $\kappa > 0$. The jump size of $P$ is in $\gamma \in [\underline{\gamma}, \overline{\gamma}]$ with constants $0 < \underline{\gamma} < \overline{\gamma} < 1$ and is generated by a probability density $g = g(\gamma)$. At each jump time $t$, the population decreases as $X_t = (1-\gamma) X_{t-}$. The harvesting rate is assumed to be a continuous-time variable having the compact range $[0, \overline{h}]$ with $\overline{h} > 0$. We assume that the decision-maker can reduce the predation pressure from $p$ to $p(1-u_t)$ with another control $u = (u_t)_{t \geq 0}$, the resource protection, valued in $[0,1]$. For example, fishery cooperatives and/or residents can control the waterbird population by executing countermeasures such as the shooting and fireworks [60, 69]. The present model therefore has the two control variables: $h$ and $u$.

The governing SDE of the population dynamics of the fish is set as

$$dX_t = -X_{t-}\left(\left(R + p(1-u_t) + h_t\right)dt + dP_t\right) \text{ for } t > 0 \text{ with } X_0 = x \geq 0. \quad (7)$$

The population dynamics are assumed to be density-independent for the sake of simplicity. The SDE (7) represents the population decrease by each of the factors explained above.

The decision-maker decides the controls $(h,u)$ to dynamically minimize

$$\phi(t, X; h, u) = \mathbb{E}^{t,x}\left[\int_t^T \left(-w_1 h_s U_s X_s + w_2 p u_s U_s X_s\right) \mathrm{d}s - w_3 U_T X_T\right], \ t > 0 \quad (8)$$

where $w_1, w_2, w_3 > 0$ are weighting factors. The first to the third terms in the right-hand side of (8) represent the cumulative utility by harvesting, the cumulative cost by taking the protection measure, and the terminal utility gained if the population remains in the river at the terminal time $t = T$. We are implicitly assuming that the harvesting cost is much smaller than the other terms. This assumption is reasonable because harvesting the fish *P. altivelis*, being different from the fisheries carried out in seas, in general uses only fishing rods [1] and/or casting nets [59] but not ships.

The last term is motivated by the fact that a larger fish is more successful in the spawning [83] and by the assumption that a more successful spawning enhances its life cycle sustainability, and thus sustainable fisheries of the fish. The liner dependence of each term on the biomass is an assumption to simplify the problem. It is more realistic to consider that the utility by the harvesting is concave with respect to $h_s U_s X_s$. Under certain conditions, each term can be replaced by concave alternatives [79] at the expense of increasing the model complexity.

The value function $\Phi = \Phi(t, x)$ as the minimized $\phi$ with respect to the controls $(h, u)$ is governed by the HJB equation

$$-\frac{\partial \Phi}{\partial t} - \inf_{h,u}\left\{\begin{array}{l}-x\left((R + p(1-u) + h)\right)\frac{\partial \Phi}{\partial x} - w_1 h U_t x + w_2 p u U_t x \\ -\kappa \int \left\{\Phi - \Phi(t, x(1-\gamma))\right\} g(\gamma) \mathrm{d}\gamma\end{array}\right\} = 0 \quad (9)$$

in $[0, T) \times (0, +\infty)$, subject to the terminal condition $\Phi_{t=T} = -w_3 U_T x$ in $(0, +\infty)$ and the boundary condition $\Phi_{x=0} = 0$ in $[0, T)$ with a polynomial growth in the far field $x \to +\infty$. The boundary and terminal conditions coincide at $(t, x) = (T, 0)$.

The HJB equation (9) admits a smooth solution of the form $\Phi(t, x) = \Psi(t) x$ with $\Psi : [0, T] \to \mathbb{R}$ if it satisfies the first-order ODE

$$-\frac{\mathrm{d}\Psi}{\mathrm{d}t} - \inf_{h,u}\left\{\begin{array}{l}-\left(\left(R + \kappa \int \gamma g(\gamma) \mathrm{d}\gamma + p(1-u) + h\right)\right)\Psi \\ -w_1 h U_t + w_2 p u U_t\end{array}\right\} = 0 \text{ in } [0, T) \quad (10)$$

subject to the terminal condition $\Psi(T) = -w_3 U_T$. It is elementary to check that, with any smooth solution $\Psi$ to (10), the function $\Psi(t) x$ solves (9). By the compactness of the ranges of the controls and the Lipschitz continuity of $U$, the ODE (10) admits a unique smooth solution. Then, $\Psi(t)x$ is continuously-differentiable with respect to both $t$ and $x$, which is therefore a smooth (and thus viscosity) solution to the HJB equation (9). As a byproduct, it is enough to estimate the net mortality rate $R + \kappa \int \gamma g(\gamma) \mathrm{d}\gamma$ but not the gradual and catastrophic ones separately. Based on this finding, in what follows, we replace $R + \kappa \int \gamma g(\gamma) \mathrm{d}\gamma$ by $R$ as a net mortality.

In summary, we can reduce an HJB equation to an ODE. We find optimal controls $(h^*, u^*)$ as functions of $t$ as follows:

$$h_t^* = \arg\min_{h \in [0, \bar{h}]} \{h(-w_1 U_t - \Psi)\} \text{ and } u_t^* = \arg\min_{u \in [0,1]} \{pu(w_2 U_t + \Psi)\}. \quad (11)$$

This control policy is advantageous from a practical viewpoint since it is free from the population, which is difficult to accurately estimate in real-world problems.

Now, the control problem reduced to solving an ODE having a smooth solution. It can be discretized with any ODE solvers such as the classical explicit and implicit Euler methods, Runge-Kutta methods, and Multi-stage methods. Here, we use the four-stage Runge-Kutta scheme. The computation here proceeds in a time-backward manner since we are considering a terminal value problem.

The parameter values are set as follows: $T = 150$ (day), $R = 0.01$ (1/day), $p = 0.01$ (1/day), $\bar{h} = 0.02$ (1/day), $w_1 = 3$, $w_2 = 2$, and different values of $w_3$. The length of the time interval is determined based on the assumption that the problem is considered from spring to autumn. For example, it corresponds to the beginning of May to the end of coming September. As a growth curve $U$, we use the classical logistic model with the growth rate 0.045 (1/day), the maximum body weight 90 (g), and the initial weight 6 (g), considering the previous study results [77, 79]. Using slightly different parameter values do not critically affect the qualitative computational results presented below. The time increment for numerical discretization of (10) is 0.01 (day), which is sufficiently smaller than $T$.

The main interest of the computation here is to see whether the real-world harvesting policy of the fish in Japan [77, 79] emerges in the presented model:

$$h_t^* = \begin{cases} 0 & (0 \leq t \leq T_0) \\ > 0 & (T_0 < t \leq T) \end{cases} \quad (12)$$

with some $T_0 \in (0, T)$. **Figures 1** and **2** show the optimal controls and the resulting population $N$, the body weight $U$, and the function $\Psi$ for $w_3 = 1$ and 3, respectively. The ranges of these variables have been normalized to the unit interval $[0,1]$ in the figure panels for the convenience of presentation. As implied in **Figure 1**, we see that there exists a range of parameter values that the control policy of the form (12), which is employed in the real world, is optimal. In this case, the countermeasure to reduce the predation pressure should be taken during the early stage. This strategy is in accordance with the fact that the mass release of the fish into a river is carried out in spring (at $t = 0$) and/or summer during which waterfowls frequently predate the fish [60].

Choosing a larger value of the parameter $w_3 = 3$ induces the optimal control to take the countermeasure both the early and late stages because of the preference of the decision-maker to make the terminal population larger. Finally, although not presented in the figure, choosing a larger value of $w_3$ such as $w_3 = 10$ leads to the optimal policy $u_t^* \equiv 1$ with $h_t^* \equiv 0$ that can maximally reduce the population decrease, but without harvesting. This is an unrealistic situation.

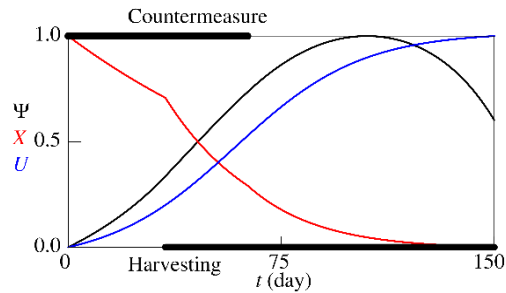

**Fig. 1**. The optimal controls controlled population $X$, the body weight $U$, and the function $\Psi$ ($w_3 = 1$). The bold lines in the upper and lower edges of the figure panel represent the time intervals during which taking countermeasure and harvesting are possible, respectively. Notice that the ranges of these variables have been normalized to the unit interval $[0,1]$.

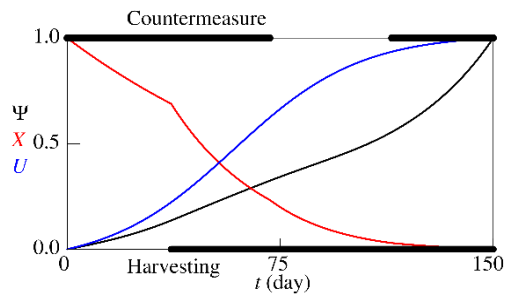

**Fig. 2**. The optimal controls controlled population $X$, the body weight $U$, and the function $\Psi$ ($w_3 = 3$). The same figure legends are the same with **Figure 1**.

### *3.2 Dam-reservoir system management*

The second problem concerns an environmentally-friendly dam-reservoir system management receiving stochastic inflows. Operating a dam-reservoir system often include multiple purposes like water resources supply, hydropower generation, and downstream environmental management, which is therefore a complex optimization problem [43]. The stochastic control approach can serve as an efficient mathematical tool for modeling and analysis of dam-reservoir systems management [65, 75, 86]. We consider a problem of a dam-reservoir system receiving stochastic inflows, as a model case to find its environmentally-friendly operation policy.

We consider the water balance dynamics of a reservoir given by

$$\mathrm{d}Y_t = (Q_t - q_t)\mathrm{d}t \text{ for } t>0 \text{ with } Y_0 = y \geq 0, \quad (13)$$

where the water volume in the reservoir at time $t$ is denoted as $Y_t$ and the outflow discharge at $t$ as $q_t$ (the control variable). The range of $Y_t$ is $\Omega = [0,\bar{Y}]$ with the prescribed volume $\bar{Y} > 0$. The inflow discharge $Q$ follows a continuous-time Markov chain $\vartheta$ with the transition matrix $\mathbf{S} = \left[ s_{i,j} \right]_{1 \leq i,j \leq I}$ with some $I \in \mathbb{N}$ [71]. The inflow discharge of the $i$ th regime is denoted as $Q(i) > 0$. Seasonality of river flows is not considered here for convenience.

The admissible range of the control $q$ must be carefully determined for rigorous mathematical modeling. Naively, the admissible range of $q$ is $A = [\underline{q},\bar{q}]$ with $0 \leq \underline{q} < \bar{q}$ imposed by technical and /or operational restrictions. We set $\underline{q} < \min_i Q(i)$ and $\bar{q} > \max_i Q(i)$ assuming a satisfactory ability of the system to handle the stochastic inflows. The range of $q$ does not have to be modulated if $Y_t \in (0,\bar{Y})$, while it has to be if $Y_t = 0, \bar{Y}$. For example, if $q_t > Q_t$ and $Y_t = 0$, then we possibly encounter the unphysical state $Y_{t+\varepsilon} < 0$ for some $\varepsilon > 0$. We set $A$ as $A = [\underline{q}, Q_t]$ ($Y_t = 0$) and $A = [Q_t, \bar{q}]$ ($Y_t = \bar{Y}$). This modification is physically appropriate, while it incurs a discontinuity of the range of controls. This difficulty can be overcome by accordingly modifying the corresponding HJB equation at the boundaries. For related problems, see Picarelli and Vargiolu [48].

The decision-maker decides the control $q$ to minimize the performance index

$$\phi(\vartheta,Y;q) = \mathbb{E}^{i,y}\left[ \int_0^\infty e^{-\delta s}\left( \frac{1}{2}(q_s - Q_s)^2 + \frac{a}{2}(\hat{q}-q_s)_+^2 + f(V_s) \right)\mathrm{d}s \right], \quad (14)$$

where $\delta > 0$ is the discount rate, $\hat{q} > 0$ is the environmental flow, and $a > 0$ is the weighting factor. The first through third terms represent penalization of the water balance condition, drawdowns from the environmental requirement with the threshold discharge $\hat{q}$ below which the downstream aquatic environment is severely affected, and a penalization of large or small water volumes, respectively. The coefficient $f \geq 0$ is Lipschitz continuous. The second term is relevant especially if the decision-maker concerns thick growth of the nuisance green filamentous algae in dam-downstream rivers due to too small river discharge [17, 35, 73].

The value function $\Phi_i = \Phi(i, y)$ as the minimized $\phi$ with respect to the control $q$ is governed by the HJB equation

$$\delta \Phi_i + \sum_{1 \leq j \leq I, j \neq i} s_{ij} \left( \Phi_i - \Phi_j \right)$$
$$-\inf_{q \in A} \left\{ \left( Q(i) - q \right) \frac{d\Phi_i}{dy} + \frac{\left( q - Q(i) \right)^2}{2} + \frac{a(\hat{q} - q)_+^2}{2} + f(v) \right\} = 0 (\geq 0) \quad \text{in } D \quad (15)$$

with $D = \{1, 2, 3, ..., I\} \times \Omega$ and the equality "=" is necessary only for $0 < y < \bar{Y}$. The controlled discharge at the boundaries $y = 0, \bar{Y}$ must satisfy the modified range as discussed above; instead, (15) is relaxed to use "$\geq$" at $y = 0, \bar{Y}$. This formulation, which is physically relevant as explained above, harmonizes with the constrained viscosity solution approach [29]. By the degenerate ellipticity of the HJB equation (15), the compactness of $A$, the signs $s_{ij} \geq 0$, and the continuity of $f$, the comparison argument (e.g., [29, Theorem 2.1]) proves the existence of at most one constrained (continuous) viscosity solution.

The local Lax-Friedrichs scheme based on the fifth-order Weighted Essentially Non-Oscillatory (WENO) reconstruction [26] combined with the fast sweeping [87] is applied to (15). This is a high-resolution scheme that has successfully been applied to degenerate elliptic equations like HJB equations. Due to the regularity deficit, meaning that solutions to the HJB equations are not smooth, the scheme seems to be less than fifth-order accurate. Nevertheless, the WENO reconstruction improves accuracy of the original local Lax-Friedrichs scheme [75].

We demonstrative a computational example assuming an existing dam-reservoir system in Japan ($\bar{Y} = 6 \times 10^7$ (m³)). The domain $\Omega$ is normalized to $[0,1]$, and $Y, q$ in the computation are normalized using $\bar{Y}$ (m³) and 86,400 (s). The parameters are the same with Yoshioka [75] where we do not consider the hydropower production here: $\underline{q} = 1.0$ (m³/s), $\bar{q} = 200$ (m³/s), $\delta = 0.01$ (1/day), $\hat{q} = 10$ (m³/s), $f = C(0.2\bar{Y} - v)_+^2 + C(v - 0.8\bar{Y})_+^2$ with $C = 5 \times 86,400/\bar{Y}$, and $a = 0.20$. The weighting factors $a, C$ have been determined to balance each term. The Markov chain $\vartheta$ employed here is the non-parametric one [75] with $I = 61$ and the discharge of $i$ th regime as $Q(i) = 1.25 + 2.5i$ (m³/s). The domain $\Omega$ is uniformly discretized with 401 vertices. The error tolerance of the fast marching is $10^{-15}$ in the $l^\infty$ norm. The initial guess is $\Phi \equiv 0$.

The numerical solution has been successfully computed using the scheme (**Figure 3**) with 8,000 iterations despite that the error tolerance of the convergence is very small. The results also present the computed optimal control $q = q^*$ normalized with respect to the inflow discharge $Q$. The decision-maker can dynamically decide the outflow discharge by referring this computational result.

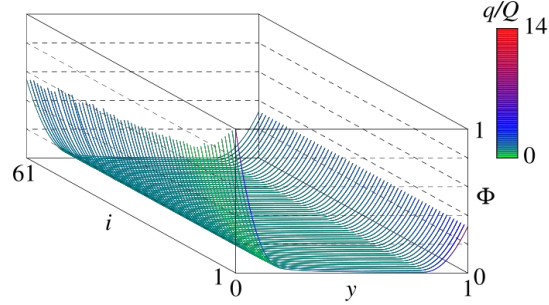

**Fig. 3**. The computed normalized value function (solid curves) and the optimal normalized discharge (color contour).

### 3.3 Algae growth management

The third problem is cost-efficient algae population management in a dam-downstream river. As explained above, thick growth of the nuisance green filamentous algae is one of the most serious environmental problems in modern river management [17, 35]. Here, we solely focus on a hydraulic control of nuisance benthic algae population dynamics following the previous study [73, 78] under a simpler setting. The model below is formulated as a deterministic control problem, but the stochastic algae population dynamics will be considered in a coupled problem later.

We consider population dynamics of benthic algae in a dam-downstream river in a lumped manner, and the population, such as the biomass per unit area of the riverbed, at time $t$ is denoted as $Z_t$. The river discharge at time $t$ as the control variable to be optimized by operating the dam is denoted as $q_t$, which has the compact range $[\underline{q}, \bar{q}]$ with $0 < \underline{q} < \bar{q}$. We assume that the population dynamics follow a control-dependent logistic model [73]:

$$\mathrm{d}Z_t = \left[ rZ_t\left(1 - \frac{Z_t}{K(q_t)}\right) - \alpha q_t Z_t \right]\mathrm{d}t, \quad t > 0, \quad Z_0 = z > 0, \tag{16}$$

where $r > 0$ is the intrinsic growth rate, $\alpha > 0$ is the proportional coefficient of detachment, and $K:[\underline{q},\overline{q}] \to \mathbb{R}_+$ is the environmental capacity as a positive, Lipschitz continuous, and increasing function. The first and second terms in the right-hand side of (16) represents the growth and the detachment by hydraulic disturbance of the population, respectively.

The unique characteristic of the model (16) is that the environmental capacity is population-dependent. As discussed in Yoshioka [73] based on the survey results [25, 64], population abundance of benthic algae and submerged vegetation can be a unimodal convex function of the flow velocity. This unimodal nature is considered because of the balance between the physical disturbance and nutrient flux transport by the river flow: too small flow discharge triggers poor algae growth while too large flow discharge triggers the algae detachment. We set the linear function $K(q_t) = K_0 + K_1 q_t$ with $K_0, K_1 > 0$ as the simplest control-dependent model. For more details, the readers should see Yoshioka [73]. With this specification, it is natural to set the range of the population as $\Omega = [0, K_0 + K_1\overline{q}]$.

The performance index to be optimized contains the disutility caused by the algae population (first term) and the penalization of the deviation between the target and realized discharge of the dam placed upstream (second term):

$$\phi(Z;q) = \mathbb{E}^z\left[\int_0^\infty \left(Z_s^m + \frac{a}{2}(q_s - \overline{q})^2\right)e^{-\delta s}\mathrm{d}s\right], \tag{17}$$

where $\delta > 0$ is the discount rate, $a > 0$ is the weighting factor, $\hat{q} \in (\underline{q},\overline{q})$ is the target discharge, and $m > 0$ is the shape parameter. The corresponding HJB equation governing the value function $\Phi = \Phi(z)$, the minimized $\phi$ by choosing $q$, is

$$\delta\Phi - \inf_{q \in [\underline{q},\overline{q}]}\left\{\left(rz\left(1 - \frac{z}{K(q)}\right) - \alpha qz\right)\frac{\mathrm{d}\Phi}{\mathrm{d}z} + \frac{a}{2}(q - \hat{q})^2 + z^m\right\} = 0 \text{ in } \Omega, \tag{18}$$

which is satisfied up to the boundary points. Notice that we can alternatively set the formal boundary condition $\Phi(0) = 0$ by taking the limit $z \to +0$ in (18). The HJB equation (18) admits at most one continuous viscosity solution by the comparison argument [73], owing to the regularity of the coefficients. We see that the uniqueness in the viscosity sense still holds true if the monomial term $Z_s^m$ in (17) is replaced by a discontinuous one such as $\chi_{\{Z_T \geq \overline{z}\}}$ with some $\overline{z} > 0$. The comparison argument in the discontinuous case can follow with a special auxiliary function in the contradiction argument (e.g., [11, Theorem 11.4]).

We proceed to numerical computation of the HJB equation (18). We use the monotone finite difference scheme combining the one-sided upwind discretization for the controlled part and the exponential discretization for the non-controlled part [73]. This scheme is monotone, stable, consistent, and is therefore convergent in a viscosity sense, implying that numerical solutions generated by the scheme converge to the unique viscosity solution (if it exists) uniformly in the (compact) domain $\Omega$ [4]. A drawback of the scheme is that it is at most first-order accurate due to the monotonicity; however, it can be combined with the policy iteration algorithm: a fast Newton-like method for solving discretized HJB equations [2]. The tri-diagonal nature of the coefficient matrix of the discretized HJB equation enhances computational efficiency by using the Thomas method [61].

The domain $\Omega$ is normalized to the unit interval $[0,1]$. The computational condition is specified below: $r=1$, $K_0=0.4$, $K_1=0.3$, $\alpha=0.5$, $\delta=2.0$, $\underline{q}=0.1$, $\overline{q}=2.0$, $\hat{q}=1.0$, $a=0.1$, and $m=0.5$. We get the maximum value of the environmental capacity $K_0+K_1\overline{q}=1$ being consistent with the length of the normalized domain defined above. The previous study assumed the convex case $m=2.0$, while the computation here assumes the concave and less regular case; the latter is expected to be concave at least near $z=0$ by invoking the asymptotic analysis result [73]. The domain $\Omega$ is divided with 501 equidistant vertices and the policy iteration is terminated if the error tolerance $10^{-14}$ in the $l^\infty$ norm is satisfied.

We compare the value function $\Phi$ (**Figure 4**) and the optimal control $q^*$ (**Figure 5**) for different values of the weighting factor $a$. Each computation terminated with less than or equal to five iterations starting from the initial guess $\Phi=0$ in $\Omega$. **Figure 4** shows that the profiles of the computed value functions $\Phi$ are concave as predicted by the asymptotic analysis result, and are difficult to distinguish with each other under the specified computational condition. **Figure 5** implies a transition of the optimal control $q^*$ between $a=0.75$ and $a=1.00$ such that $q^*$ seems to be continuous and close to the target value $\hat{q}=1$ for relatively large $a\geq 1$, while it is discontinuous and clearly different from $\hat{q}=1$ for smaller $a$. This kind of sudden transition is due to that the minimizer of the second term of (18), assuming that it is an interior solution ($q^*\in(\underline{q},\overline{q})$), is a solution to a non-monotone third-order polynomial. This kind of phenomenon does not occur if $K_0=1$ and $K_1=0$, as demonstrated in **Figure 6**. The non-trivial profile of the optimal control $q^*$ in the present control-dependent model is due to considering the balance between the physical disturbance and nutrient transport, which is not considered in the standard control-independent model.

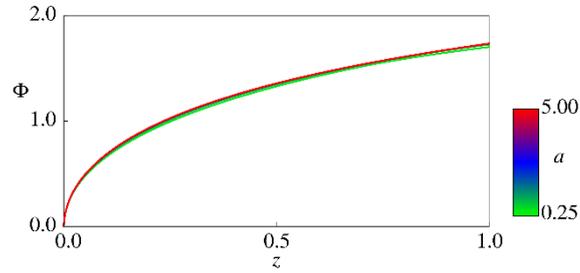

**Fig. 4.** The computed value functions $\Phi$ for the weighting factors $a = 0.25i$ ($i = 1, 2, 3, ..., 20$).

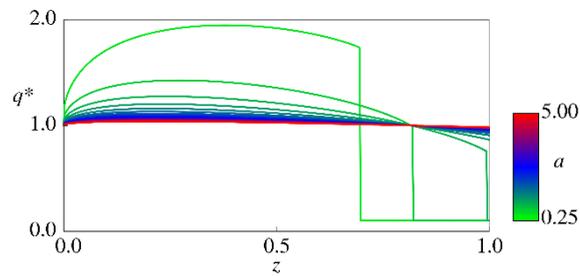

**Fig. 5.** The computed optimal control $q^*$ for the weighting factors $a = 0.25i$ ($i = 1, 2, 3, ..., 20$).

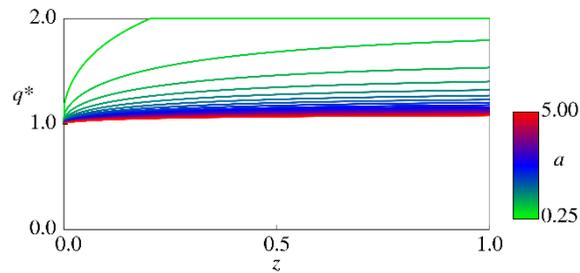

**Fig. 6.** The computed optimal control $q^*$ for the weighting factors $a = 0.25i$ ($i = 1, 2, 3, ..., 20$).with the control-independent environmental capacity with $K_0 = 1$ and $K_1 = 0$.

### 3.4 Sediment storage management

The last example is a sediment replenishment problem in a dam-downstream river. This example is different from the previous ones in the sense that the decision-maker can intervene only discretely and randomly, whereas he/she can control the target dynamics continuously in the previous problems. It is still not always possible to collect information of environmental and ecological dynamics under natural environment and are often provided only discretely like weekly or monthly [55, 68].

The sediment replenishment problem in a dam-downstream river considered here has a simple environmental background; sediment supply, which usually occurs along a natural river, is stopped at a dam. A critical issue is that the benthic community in the dam-downstream river is critically affected by the absence of the sediment supply [19, 41, 54]. Some environmental managers encountering this issue have been trying to replenish earth and soils from outside the river [57]. However, cost-effectiveness of the sediment replenishment has only recently been considered from a mathematical side by the authors [74, 81-82, 84].

Assume that we can place a sediment lump in a dam-downstream reach. The storable amount of the sediment is $\bar{W} > 0$. The amount of stored sentiment at time $t$ is denoted as $W_t$ with the range $\Omega = [0, \bar{W}]$. Assume a constant river flow having a sufficiently large discharge such that the sediment is transported toward downstream. The transport rate as a function of the hydraulic variables is then considered as a constant $S > 0$ [67, 74]. The corresponding sediment storage dynamics are

$$\mathrm{d}W_t = -S\chi_{\{W_t > 0\}}\mathrm{d}t \text{ for } t > 0 \text{ with } W_0 = w \in \Omega. \tag{19}$$

This is a discontinuous dynamical system [15] where the drift coefficient is discontinuous when $W_t = 0$. It has a unique continuous solution in the Filippov sense: $W_t = \max\{0, W_0 - St\}$ ($t \geq 0$) owing to the one-sided Lipschitz property of (19).

We assume that the decision-maker can replenish the sediment storage to the maximum value impulsively. This assumption is valid if the sediment replenishment can be carried out with a much shorter time than the time-scale of the dynamics. We further assume that the intervention can be carried out only discretely and randomly, and the chances of the interventions are identified as the jump times of a Poisson process $N = (N_t)_{t \geq 0}$ with the intensity $\lambda > 0$. The inverse $\lambda^{-1}$ then gives the mean as well as standard deviation of the intervals between successive interventions. The quantity $\lambda^{-1}$ is therefore the characteristic time scale of interventions.

The amount of sediment replenished at time $t$ is denoted as $\eta = (\eta_t)_{t \geq 0}$ (For its deeper mathematical details, see [63]). Clearly, the state that the decision-maker should avoid is the sediment depletion $W_t = 0$. In addition, sediment replenishment can be costly. We assume that the decision-maker must pay both fixed cost $d > 0$ (e.g., labor cost) and the proportional cost $c\eta_t$ when he/she replenishes sediment.

Based on the consideration above, the performance index is set as

$$\phi(W, \eta) = \mathbb{E}^w\left[\int_0^\infty \chi_{\{W_s = 0\}}e^{-\delta s}\mathrm{d}s + \int_0^\infty (c\eta_s + d)\chi_{\{\eta_s > 0\}}e^{-\delta s}\mathrm{d}N_s\right] \tag{20}$$

with $\delta > 0$, where the first term penalizes the sediment depletion, while the second term is the incurred costs. By a dynamic programming argument [63], the value function $\Phi = \Phi(w)$ as the minimized $\phi$ with respect to $\eta$ is governed by

$$\delta\Phi + S\chi_{\{w>0\}}\frac{d\Phi}{dw} + \lambda\left(\Phi - \min_{\eta\in\{0,1-w\}}\left\{\chi_{\{\eta>0\}}(c\eta+d) + \Phi(w+\eta)\right\}\right) = \chi_{\{w=0\}} \text{ in } \Omega. \quad (21)$$

No boundary condition is necessary since the HJB equation is assumed to be satisfied up to the boundary points. The optimal replenishment policy is

$$\eta_t^* = \eta^*(W_t) = \arg\min_{\eta\in\{0,1-w\}}\left\{\chi_{\{\eta>0\}}(c\eta+d) + \Phi(W_t+\eta)\right\}. \quad (22)$$

This applies only at the randomly arriving intervention chances.

Despite that the coefficients in (21) are non-smooth, Yoshioka [74] proved that the HJB equation admits a smooth exact solution and that it is the value function under certain condition of the parameters. The optimal replenishment policy is

$$\eta^*(W_t) = (1-W_t)\chi_{\{W_t\leq\bar{w}\}}, \quad (23)$$

representing a strategy that the sediment should be replenished if its storage is below the threshold $\bar{w}$. Finding the optimal control thus reduces to finding $\bar{w}$ if it exists.

We present numerical examples with a regularized counterpart where the characteristic function $\chi_{\{w=0\}}$ is replaced by $\chi_\varepsilon(w) = \max\{0, 1-w/\varepsilon\}$ with a regularization parameter $0 < \varepsilon < 1$. Similarly, we replace $\chi_{\{w>0\}}$ by the regularized one $1 - \chi_\varepsilon(w)$. This regularization, which has not been used so far, is introduced for the two reasons. The first reason is to guarantee the unique existence of continuous viscosity solutions. The verification argument [74] implies that a smooth solution is a value function; however, the unique existence of (viscosity) solutions was not discussed. The regularization method does not critically affect the dynamics and optimal control, as demonstrated below, while it guarantees that the equation (21) has Lipschitz continuous coefficients if $\varepsilon > 0$, and we see that the comparison theorem applies [11, Chapter 3] with a slight adaptation to the presented model. The second reason is to improve convergence of numerical schemes expecting that viscosity solutions have higher regularity with smoother coefficients.

We present a computational example of the HJB equation (21) with the regularization method to numerically show that the optimal policy is the threshold type and the impacts of the regularization. The parameter values are set as follows: $S = 0.1$, $\lambda = 0.1$, $c = 0.5$, $d = 0.4$, and $\delta = 0.1$. The HJB equation (21) is numerically dis-cretized using the semi-Lagrangian scheme with the third-order WENO interpolation [13]. Its advantage is high stability owing to the semi-Lagrangian nature combined with the high-order and sharp interpolation by the WENO reconstruction. Its disadvantage is that a careful scaling required between the temporal increment $\Delta t$ and spatial increment $\Delta x$ for convergence [13]. In our case, the scheme is 1.5th-order accurate if $\Delta t$ is proportional to $(\Delta x)^{1.5}$. We use $\Delta x = 1/300$ and $\Delta t = 30(\Delta x)^{1.5}$. The computation starts from the initial guess $\Phi \equiv 0$ and employs a value function iteration with the error tolerance of $10^{-10}$ in the $l^\infty$ norm.

We examine the regularization parameters $\varepsilon = 0.1, 0.05, 0.01, 0.005,$ and $0.001$. **Figures 7** and **8** show the computed value functions $\Phi$ and the auxiliary variables $\omega$: $\omega = 1$ if $\eta^* = 1 - w$ and $\omega = 0$ otherwise. The results with $\varepsilon = 0.005$ and $0.001$ are very close to that with $\varepsilon = 0.01$, and are therefore not plotted in the figures. Regularizing the non-smooth coefficients locally influence the solution shape near $x = 0$ and globally affects the magnitude in the domain. The impacts of regularization get smaller as $\varepsilon$ decreases. It should also be noted that **Figure 8** implies that the optimal control is indeed the threshold type under the regularization. The computational results suggest choosing $\varepsilon = O(\Delta x)$ to compute reasonably accurate numerical solutions. We expect the convergence of numerical solutions to the unique viscosity solution under this limit; this is an open issue.

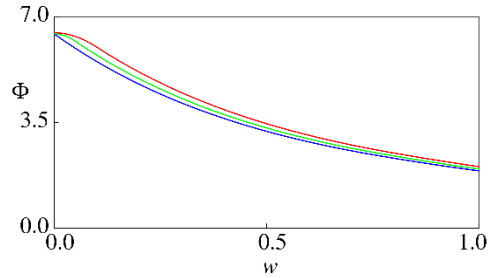

**Fig. 7.** The computed value functions $\Phi$ for $\varepsilon = 0.1$ (Red), 0.05 (Green), and 0.01 (Blue).

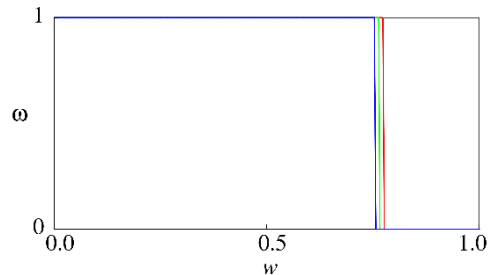

**Fig. 8.** The computed auxiliary functions $\omega$ for $\varepsilon = 0.1$ (Red), 0.05 (Green), and 0.01 (Blue).

## 4. Coupled problem

*4.1 Overview*

We discussed only 1-D examples above. In this section, we present a demonstrative example of formulation and computation of a coupled fisheries and environmental problem. This is a control problem of a dam-reservoir system and its downstream environment where the water balance dynamics receiving regime-switching stochastic inflows, nuisance benthic algae population dynamics, and lumped sediment storage dynamics are concurrently considered. The fishery resource dynamics are considered in a performance index.

Due to high-dimensional nature of the problem having the four state variables in addition to the time variable, numerical computation of its associated HJB equation is very costly if we rely on a standard numerical method like a finite difference scheme and a semi-Lagrangian scheme on a conventional structured grid. More specifically, such conventional numerical methods encounter exponential increase of the computational cost with respect to the increase of the increase of the total number of state variables: the curse-of-dimensionality. We employ the sparse grid technique to alleviate this computational issue [10].

*4.2 Control problem*

We present a management problem concerning a dam-reservoir system. This is a control problem based on the above-discussed separate problems. As explained above, regulated river flows downstream of a dam encounter a critical reduction of the sediment transport and associated environmental problems like the thick growth of nuisance benthic algae. The latter potentially leads to the reduction of the growth rate of fishery resources like *P. altivelis*. From a fisheries viewpoint, mitigating the environmental impacts on the growth of the fish can be addressed at least by controlling the outflow discharge of the dam and/or replenishing the sediment.

The system dynamics to be controlled have the four state variables: the inflow of the reservoir as a continuous-time Markov chain $\vartheta$ and the three continuous-time variables, which are the water volume of the reservoir $X_1$, the sediment storage $X_2$, and the nuisance benthic algae population $X_3$. The total number of the regimes of $\alpha$ is denoted as $I \in \mathbb{N}$ and the range of the variable $X_j$ as $\Omega_j = [0, \bar{X}_j]$ ($j = 1, 2, 3$). The governing system of $X_j$ ($j = 1, 2, 3$) as follows:

$$dX_{1,t} = (Q_t - q_t)dt$$
$$dX_{2,t} = -S(q_t)\chi_{\{X_{2,t}>0\}}dt + \eta_t dN_t \quad \text{for } t > 0 \quad (24)$$
$$dX_{3,t} = \left[rX_{3,t}\left(1 - \frac{X_{3,t}}{K(q_t)}\right) - \alpha(X_{2,t})q_t X_{3,t}\right]dt$$

with initial conditions $X_{j,0} \in \Omega_j$ ($j = 1, 2, 3$) and $\alpha_0 \in i$, where the same notations of the parameters with the previous section are utilized except for the transport rate $S$ and detachment coefficient $\alpha$. The transport rate $S$ was assumed to be a constant in the previous section but is now a non-negative coefficient depending on the outflow discharge. The detachment coefficient $\alpha$ was also set as a constant, but now as a function $\alpha(X_{2,t})$ depending on the sediment storage. We assume $\alpha(0) = 0$ because the algae detachment is not significant if the river flow does not contain soil particles (bedload) [84, Chapter 3]. Set the sequence representing the regimes as $\Theta = \{i\}_{1 \leq i \leq I}$.

The control variables of the present problem are the outflow discharge $q = (q_t)_{t \geq 0}$ as in the second model problem and the sediment replenishment $\eta = (\eta_t)_{t \geq 0}$ as in the fourth problem. We assume that, at each time $t$, the outflow discharge is valued in the discrete set $A_t = \{a_j Q_t\}_{0 \leq j \leq n}$ with some $n \in \mathbb{N}$ and some real $a_j$. The range $A_t$ is modified appropriately if $X_{t,1} = 0$ or $X_{t,1} = \overline{X}_1$ as discussed in **Section 3.2**. The reservoir volume $X_1$ is then a.s. confined in $\Omega_1$. As in **Section 3.4**, the replenishment is specified so that it equals either 0 (Do nothing) or $\overline{X}_2 - X_{2,t}$ (Fully replenish, only when $\overline{X}_2 < X_{2,t}$) at each jump time $t$ of the Poisson process $N$ representing the a sequence of chances to replenish the sediment.

We consider a performance index containing a penalization of the state variables deviating from a compact safe region $\Omega_S \subset \Omega_1 \times \Omega_2 \times \Omega_3$ (first term) at a terminal time $T > 0$, a penalization of the deviation of the outflow discharge from the inflow (second term), and a cumulative sediment replenishment cost (third term):

$$\phi(t, \vartheta, X_1, X_2, X_3; \eta, q)$$
$$= \mathbb{E}^{t,i,x_1,x_2,x_3}\left[\begin{array}{l} \int_t^T e^{-\delta(s-t)} f(X_{1,s}, X_{2,s}, X_{3,s})ds \\ \int_t^T e^{-\delta(s-t)} \frac{1}{2}\left(1 - \frac{q_s}{Q_s}\right)^2 ds + \int_t^T (c\eta_t + d)\chi_{\{\eta_t > 0\}} e^{-\delta(s-t)} dN_s \end{array}\right], \quad (25)$$

with the discount rate $\delta \geq 0$, coefficients of proportional and fixed costs $c, d > 0$, and a continuous function $f$ such that it equals 0 if its arguments fall in the safe region $\Omega_S$ and positive otherwise. This performance index means that the river environmental condition should be in the safe region by controlling the system dynamics. The interval $[0, T]$ can be chosen as the time around which the fish *P. altivelis* in the dam-downstream river stars to significantly growing up. At that time, the algae population, sediment storage, and reservoir water volume, should be in some safe region; otherwise, the growth of the fish can be critically affected. However, adjusting the river environmental condition to the required state by operating the dam-reservoir system and/or replenishing the sediment is costly. The problem considered here is only one example of the coupled modeling, but many other situations can be covered by modifying the performance index.

The value function is the minimized $\phi$ with respect to the controls $q$ and $\eta$:

$$\Phi(t, i, x_1, x_2, x_3) = \Phi_i(t, x_1, x_2, x_3) = \inf_{q, \eta} \phi(t, \vartheta, X_1, X_2, X_3; \eta, q). \tag{26}$$

The associated HJB equation that governs $\Phi$ is

$$-\frac{\partial \Phi_i}{\partial t} + \delta \Phi_i + \sum_{1 \leq j \leq I, j \neq i} s_{ij}(\Phi_i - \Phi_j) - f(x_1, x_2, x_3)$$

$$-\inf_q \left\{ \begin{array}{l} \left(Q(i) - q\right)\frac{\partial \Phi}{\partial x_1} - S(q)\chi_{\{x_2 > 0\}}\frac{\partial \Phi}{\partial x_2} + \left(rx_3\left(1 - \frac{x_3}{K(q)}\right) - \alpha(x_2)qx_3\right)\frac{\partial \Phi}{\partial x_3} \\ + \frac{1}{2}\left(1 - \frac{q}{Q(i)}\right)^2 \end{array} \right\} \tag{27}$$

$$+ \lambda\left(\Phi - \min_{\eta \in \{0, \bar{X}_2 - x_2\}}\left\{\chi_{\{\eta > 0\}}(c\eta + d) + \Phi(t, i, x_1, x_2 + \eta, x_3)\right\}\right) = 0$$

in $(0, T) \times \Theta \times \Omega_1 \times \Omega_2 \times \Omega_3$, subject to the terminal condition

$$\Phi(T, i, x_1, x_2, x_3) = 0 \text{ in } \Theta \times \Omega_1 \times \Omega_2 \times \Omega_3. \tag{28}$$

It seems to be hopeless to analytically handle the minimization term (27) because of its complexity, implying that we have to tackle this issue numerically.

### *4.3 Numerical scheme*

We use the semi-Lagrangian scheme on sparse grids [8]. The sparse grid technique does not utilize the conventional fully tensorized computational grids but only grids containing selected basis functions to balance computational costs and accuracy. Several sparse grids are available depending on the problem to be considered [10, 28, 53]. Computational complexity of spars grids is quasi-optimal in terms of theoretical computational efficiency with respect to $L^2$- and $L^\infty$-norms, and are more suited to approximating higher-dimensional functions[10, Section 3].

We employ the modified sparse grids [28] that are comparably accurate with the standard ones [10] but have a relatively smaller number of the grid points on boundary. This means that the former exhibit higher accuracy in the interior of the domain if the function to be approximated are sufficiently smooth near and along the boundary. The implementation method and a series of theoretical analysis results on this sparse grid technique are available in Kang and Wilcox [28]. Later, **Figure 9** shows that the sparse grid we use is indeed different from the standard grids.

An important point of the sparse grid technique is its flexibility to approximate value functions arising in high-dimensional problems; our HJB equation is not an exception. Another important point is that its approximation error in the $L^2$-sense is the order of $N^{-2}(\log N)^{d-1}$ that is slightly worse than the order $N^{-2}$ of fully-tensorized (conventional) grids, where $N = 2^{l-d}$ is the degree-of-freedom per dimension, $l \in \mathbb{N}$ ($l \geq d$) is discretization level of the sparse grids, and $d$ is the total number of dimensions. The grid becomes finer as the discretization level $l$ increases. The implementation cost of the sparse grids is scaled as $N(\log N)^d$, while that of the standard grids as $N^d$, implying that the total cost of the former is far smaller than that of the latter in return for the slight accuracy deterioration.

Functions approximated on sparse grids should be at least twice partially differentiable in the weak sense [10, Theorem 3.8]. Too much irregular functions, like that having a discontinuity, should not be handled directly on sparse grids. In addition, the sparse grid technique does not guarantee monotonicity of the interpolations on them, meaning that some limitation method may be necessary in applications especially if the function to be approximated are not smooth [66]. We do not use limitation methods here but will be investigated in future.

### 4.4 Computational conditions

The coefficient and parameter values are specified considering the approximation ability of the sparse grids. There are the four coefficients to be specified: $K$, $S$, $\alpha$, and $f$. For the sake of simplicity, we set the constant environmental capacity case $K = 1$ with which the fastest logistic growth of the nuisance algae occurs. In addition, we normalize each $\Omega_j$ to the unit interval [0,1] by appropriately normalizing each state variable. The Markov chain $\vartheta$ utilized here is based on the same data and identification method but has the smaller total number of regimes with that in **Section 3.2**. Here, we set $I = 21$ and $Q(i) = 2.5 + 5(i-1)$ (m$^3$/s).

The transport rate $S$ is based on the semi-empirical physical formulae based on the [67, 74] with the following specified hydraulic condition: river width 30 (m), longitudinal riverbed slope 0.001, soil (sediment) particles diameter 0.006 (m), Manning's roughness coefficient 0.03 (m$^{-1/3}$/s), the density of water 1,000 (kg/m$^3$), the soil density 2,600 (kg/m$^3$), and the gravitational acceleration coefficient 9.81 (m/s$^2$). Using these values, the transport rate $S$ (m$^3$/day) is calculated as a function of the outflow discharge $q$ (m$^3$/s) under the normalization of the state variables as

$$S(q) = \bar{X}_3^{-1} \times A \times \max\{Bq^{0.6} - C, 0\}^{1.5} \quad (1/\text{day}) \tag{29}$$

with the constants $A = 3.82 \times 10^4$ (m$^3$/s), $B = 1.31 \times 10^{-2}$ (s/m$^3$), and $C = 4.7 \times 10^{-2}$. The maximum volume of the storable sediment in the dam-downstream river reach, which is $\bar{X}_3$, is assumed to be 200 (m$^3$).

The coefficient $\alpha$ of the algae detachment is a function of the sediment storage $x_2$ such that the detachment occur if the sediment is not depleted ($x_2 > 0$). In addition, several field survey results suggest that the channel bed disturbance by the bedload transport is a driver of the algae detachment [30, 36]. As a simple model, we propose to set $\alpha(x_2) = \alpha_0 x_2^m$ (s/m$^3$/day) with positive constants $\alpha_0$ and $m$. We set $\alpha_0 = 0.1$ (s/m$^3$/day) and $m = 0.5$. The growth rate $r$ is set as 0.5 (1/day) considering the previous research results [73, 78].

The coefficient $f$ of the terminal condition must be carefully specified because of the approximation ability of the sparse grids. We consider $f$ of the form

$$f = \begin{bmatrix} \max\{x_1 - \bar{x}_1, 0\}^p / (1 - \bar{x}_1)^p + \max\{\underline{x}_1 - x_1, 0\}^p / \underline{x}_1^p \\ + \max\{\underline{x}_2 - x_2, 0\}^p / \underline{x}_2^p + \max\{x_3 - \bar{x}_3, 0\}^p / (1 - \bar{x}_3)^p \end{bmatrix} \quad (1/\text{day}) \tag{30}$$

with positive constants $\bar{x}_1, \underline{x}_1, \underline{x}_2, \bar{x}_3$, and $p$. This $f$ means that the reservoir volume should be neither too small nor too large, the algae population should not be large, and the sediment storage should not be small. We choose $p = 3$ with which the smoothness $f \in C^2(\Omega_1 \times \Omega_2 \times \Omega_3)$ holds true and therefore this $f$ can be approximated on the sparse grids [10, Theorem 3.8]. We set $\bar{x}_1 = 0.8$, $\underline{x}_1 = 0.2$, $\underline{x}_2 = 0.2$, $\bar{x}_3 = 0.8$. The coefficients $\{a_j\}_{0 \leq j \leq n}$ of the control set $A_t$ is set as $a_0 = 0$, $a_1 = 1.0/2.0$, $a_2 = 2.0$, $a_3 = 1.0/3.0$, and $a_4 = 3.0$. The coefficients on the replenishment are set as $c = 0.15$ and $d = 0.05$.

The time horizon is set as $T = 60$ (day) assuming a management problem in a growth period (late spring to the coming summer in a year) of the fish *P. altivelis*. The time increment for the temporal integration of the semi-Lagrangian scheme is set as 0.005 (day), and we used the modified sparse grid of the level 11 in the sense of Kang and Wilcox [28] for discretization of the 3-D space $\Omega_1 \times \Omega_2 \times \Omega_3$. The total number of grid points is 6,017, and the minimum distance among the grid points is 1/256 (**Figure 9**). The corresponding full-tensor grid requires $O(10^7)$ grid points, which is greater than that utilized here. Furthermore, we have 21 regimes, essentially implying that the total numbers of grid points are $O(10^5)$ and $O(10^8)$; the former is about $O(10^3)$ times smaller than the latter, demonstrating high efficiency of the sparse grid technique.

### 4.5 Numerical computation

The obtained numerical solution at several time steps for small and large inflow discharges are presented to analyze the optimal controls as functions of the state variables. **Figures 10** through **12** show the computed value functions $\Phi$, the optimal sediment replenishment ($\eta^*$), and the optimal level of the outflow discharge ($q^*$) of the relatively low and high flow regimes at several time instances. The computational results suggest that replenishing the sediment in this case is in general optimal if the sediment storage is not the full. A large water storage in the reservoir or a small alga population under a relatively high flow regime is found to be another possible condition where replenishing the sediment is not optimal.

The optimal level of the outflow discharge highly depends on the state variables, and it seems to be not easy to find a simple law governing it. Nevertheless, at least theoretically, the decision-maker can decide the control variables at each time instance by referring to the computational results over the phase space. A technical issue is that visualizing a high-dimensional data is usually a difficult task. In our example, each figure panel of **Figures 10** through **12** is only the plot at some instance and some regime. Establishment of an effective visualization technique for high-dimensional data would be an interesting research topic related to optimization and control of many engineering problems. Utilizing some explainable artificial intelligence technique may be beneficial for better understanding the high-dimensional data like the output of this stochastic control model.

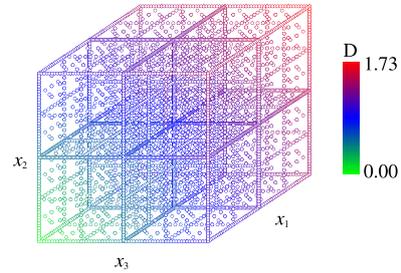

**Fig. 9.** The sparse grid. The color represents the distance ("D" in the figure) from the origin.

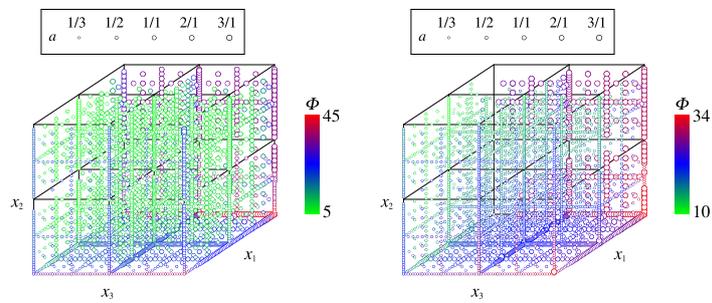

**Fig. 10.** The computed optimal controls at the time $t = 0$ (day) for the relatively low (Left: $i = 3$) and high (Right: $i = 13$) inflow regimes. Only the points where the sediment should be replenished ($\eta^* > 0$) are plotted.

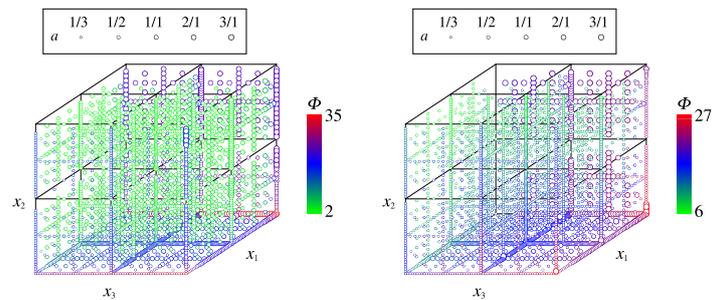

**Fig. 11.** The computed optimal controls at the time $t = 30$ (day) for the relatively low (Left: $i = 3$) and high (Right: $i = 13$) inflow regimes. The same legends with **Figure 10**.

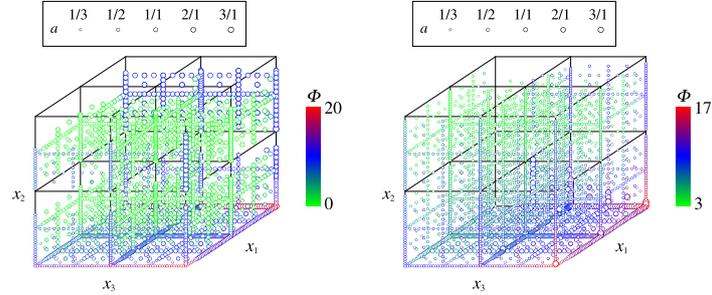

**Fig. 12.** The computed optimal controls at the time (day) for the relatively low (Left: $i = 3$) and high (Right: $i = 13$) inflow regimes.

## 5. Conclusions

We analyzed independent stochastic control problems on fishery and environmental problems from a unified viewpoint. All the problems based on a dynamic programming argument, which reduces solving a stochastic control problem to solving some HJB equation. We discussed that these problems can be analyzed from the viewpoint of continuous viscosity solutions and demonstrated computational examples. Finally, we considered a coupled problem utilizing a semi-Lagrangian scheme with the sparse grid technique.

The problems we handled are only the tip of the iceberg: many other interesting topics remain to be addressed. For example, a stochastic control problem with more than one decision-makers would arise when considering integrated management of a watershed where many industries, such as agriculture, fisheries, urban planning, coexist. A massive numerical computation architecture as well as a smart numerical technique like the sparse grid technique should be combined to tackle this advanced issue. Assessing viability of system control strategies is also an important topic in engineering applications [21]. Viability analysis of high-dimensional deterministic dynamical systems have been analyzed [9], but the analysis of stochastic dynamical systems seems to be far less analyzed. We are currently tackling this issue by using the proposed framework based on the dynamic programming methodology and numerical schemes. Employing a multi-objective formulation [18] would facilitate modeling and analysis of the stochastic viability.

Evaluating economic impacts of inland fisheries through river environmental management would be essential for achieving Sustainable Development Goals (SDGs) and symbiosis of humans and environment. A key would be setting a utility function to be optimized from both mathematical and engineering viewpoints. A collaboration of science and engineering will be necessary to resolve these complex but emergent issues. We hope that the presented mathematical framework would become a cornerstone to approach these issues.

**Acknowledgments** JSPS Research Grant No. 19H03073 and Kurita Water and Environment Foundation Grant No. 19B018 and No. 20K004, and grants from MLIT Japan for ecological survey of a life history of the landlocked *P. altivelis* and management of seaweed in Lake Shinji support this research. The author gratefully thanks to Dr. Srikanta Patnaik for his invitation to this volume.